\documentclass[11pt]{amsart}
\usepackage{amssymb,amsmath,graphics,graphicx,amscd}

\newtheorem{theorem}{Theorem}[section]
\newtheorem{lem}[theorem]{Lemma}
\newtheorem{prop}[theorem]{Proposition}
\newtheorem{cor}[theorem]{Corollary}

\theoremstyle{definition}

\theoremstyle{remark}

\numberwithin{equation}{section}

\newcommand{\lam}{\lambda}

\newcommand{\Hor}{{\mathcal{H}}}

\newcommand{\ra}{\rightarrow}

\newcommand{\Ad}{\text{Ad}}

\newcommand{\lb}{\langle}
\newcommand{\rb}{\rangle}

\newcommand{\mg}{\mathfrak{g}}
\newcommand{\mh}{\mathfrak{h}}

\renewcommand{\mp}{\mathfrak{p}}
\newcommand{\bs}{\backslash}

\newcommand{\R}{\mathbb{R}}

\begin{document}

\newcommand{\spacing}[1]{\renewcommand{\baselinestretch}{#1}\large\normalsize}
\spacing{1.14}

\title{Invariant metrics with nonnegative curvature on compact Lie Groups}

\author{Nathan Brown}
\email{nbrown@northwestern.edu}
\author{Rachel Finck}
\email{rfinck@mit.edu}
\author{Matthew Spencer}
\email{matthew.p.spencer@williams.edu}
\author{Kristopher Tapp}
\email{ktapp@williams.edu}
\author{Zhongtao Wu}
\email{ztwu@mit.edu}
\thanks{Supported in part by NSF grant DMS--0303326 and NSF grant DMS-0353634.}

\date{\today}

\begin{abstract}
We classify the left-invariant metrics with nonnegative sectional curvature on $SO(3)$ and $U(2)$. 
\end{abstract}

\maketitle


\section{Introduction}\label{intro}
Every compact Lie group admits a bi-invariant metric, which has nonnegative sectional curvature.
This observation is the starting point of almost all known constructions of nonnegatively and positively curved manifolds.  The only exceptions are found in~\cite{Cheeger}, \cite{GZ}, and \cite{Yang}, and many of these exceptions admit a different metric which does come from a bi-invariant metric on a Lie group; see~\cite{totaro}.

In order to generalize this cornerstone starting point, we consider the following problem: for each compact Lie group $G$, classify the left-invariant metrics on $G$ which have nonnegative sectional curvature.  The problem is well-motivated because each such metric, $g_l$, provides a new starting point for several known constructions.

For example, suppose $H\subset G$ is a closed subgroup.  The left-action of $H$ on $(G,g_l)$ is by isometries.  The quotient space $H\bs G$ inherits a metric of nonnegative curvature, which is generally inhomogeneous.  Geroch proved in~\cite{ger} that this quotient metric can only have positive curvature if $H\bs G$ admits a normal homogeneous metric with positive curvature, so new examples of positive curvature cannot be found among such metrics.  However, it seems worthwhile to search for new examples with quasi- or almost-positive curvature (meaning positive curvature on an open set, or on an open and dense set).

A second type of example occurs when $H\subset G$ be a closed subgroup, and $F$ is a compact Riemannian manifold of nonnegative curvature on which $H$ acts isometrically.  Then $H$ acts diagonally on $(G,g_l)\times F$.  The quotient, $M=H\bs(G\times F)$, inherits a metric with nonnegative curvature.  $M$ is the total space of a homogeneous $F$-bundle over $H\bs G$.  Examples of this type with quasi-positive curvature are given in~\cite{Tapp}.  It is not known whether new positively curved examples of this type could exist.

In this paper, we explore the problem for small-dimensional Lie groups.  In section 3, we classify the left-invariant metrics on $SO(3)$ with nonnegative curvature.  This case is very special because the bi-invariant metric has positive curvature, and therefore so does any nearby metric.  In section 4, we classify the left-invariant metrics on $U(2)$ with nonnegative curvature.  This classification does not quite reduce to the previous, since not all such metrics lift to product metrics on $U(2)$'s double-cover $S^1\times SU(2)$.

The only examples in the literature of left-invariant metrics with nonnegative curvature on compact Lie groups are obtained from a bi-invariant metric by shrinking along a subalgebra (or a sequence of nested subalgebras).  We review this construction in section 5, which in its greatest generality is due to Cheeger~\cite{Cheeger}.  Do all examples arise in this way?  In section 6, we outline a few ways to make this question precise.  Even for $U(2)$, the answer is essentially no. Also, since each of our metrics on $U(2)$ induces via Cheeger's method a left-invariant metric with nonnegative curvature on any compact Lie group which contains a subgroup isomorphic to $U(2)$, our classification for small-dimensional groups informs the question for large-dimensional groups.

The fourth author is pleased to thank Wolfgang Ziller for helpful suggestions.
\section{Background}
Let $G$ be a compact Lie group with Lie algebra $\mg$.  Let $g_0$ be a bi-invariant metric on $G$.  Let $g_l$ be a left-invariant metric on $G$.  The values of $g_0$ and $g_l$ at the identity are inner products on $\mg$ which we denote as $\lb\cdot,\cdot\rb_0$ and $\lb\cdot,\cdot\rb_l$.  We can always express the latter in terms of the former:
$$\lb A,B\rb_l=\lb\phi(A),B\rb_0$$
for some positive-definite self-adjoint $\phi:\mg\ra\mg$.  The eigenvalues and eigenvectors of $\phi$ are called eigenvalues and eigenvectors vectors of the metric $g_l$.

There are several existing formulas for the curvature tensor of $g_l$ at the identity, which we denote as $R:\mg\times\mg\times\mg\ra\mg$.  P\"uttmann's formula from~\cite{put} says that for all $x,y,z,w\in\mg$,

\begin{align}\label{put}
\lb R(x,y)z,w\rb_l=&-\frac{1}{4}\left(\langle\lbrack\phi x,y\rbrack,\lbrack z,w\rbrack\rangle_0+\langle\lbrack x,\phi y\rbrack,\lbrack z,w\rbrack\rangle_0\right.\\
&\text{\hspace{.34in}}\left.+\langle\lbrack x,y\rbrack,\lbrack\phi z,w\rbrack\rangle_0+\langle\lbrack x,y\rbrack,\lbrack z,\phi w\rbrack\rangle_0\right)\notag\\
&-\frac{1}{4}\left(\langle\lbrack x,w\rbrack,\lbrack y,z\rbrack\rangle_l-\langle\lbrack x,z\rbrack,\lbrack y,w\rbrack\rangle_l-2\langle\lbrack x,y\rbrack,\lbrack z,w\rbrack\rangle_l\right)\notag\\
&-\left(\langle B(x,w),\phi^{-1}B(y,z)\rangle_0-\langle B(x,z),\phi^{-1}B(y,w)\rangle_0\right),\notag
\end{align}
where $B$ is defined by $B(x,y)=\frac{1}{2}\left(\lbrack x,\phi y\rbrack+\lbrack y,\phi x\rbrack\right).$

Sectional curvatures can be calculated from this formula.  An alternative formula by Milnor for the sectional curvatures of $g_l$ is found in~\cite{milnor}.  If $\{e_1,...,e_n\}\subset\mg$ is a $g_l$-orthonormal basis, then the sectional curvature of the plane spanned by $e_1$ and $e_2$ is:
\begin{align}\label{miln}
\kappa(e_1,e_2) =\sum_{k=1}^n & \left(   \frac{1}{2}\alpha_{12k}\left(-\alpha_{12k}+\alpha_{2k1} +\alpha_{k12}\right)\right.\\
  & \left.-\frac{1}{4}\left(\alpha_{12k}-\alpha_{2k1}+\alpha_{k12}\right)\left(\alpha_{12k} +\alpha_{2k1}-\alpha_{k12}\right)-\alpha_{k11}\alpha_{k22}\right),\notag
\end{align}
where $\alpha_{ijk}=\lb[e_i,e_j],e_k\rb_l$ are the \emph{structure constants}.

\section{Metrics on $SO(3)$}
Let $G=SO(3)$, which has Lie algebra $\mg=so(3)$.  We use the following bi-invariant metric:
$$\lb A,B\rb_0 = (1/2)\text{trace}(AB^{T}) \text{ for }A,B\in\mg.$$
A left-invariant metric $g_l$ on $G$ has three eigenvalues, $\{\lambda_1^2,\lambda_2^2,\lambda_3^2\}$.  The metric can be re-scaled so that $\lambda_3=1$. 
\begin{prop}\label{P:rocket}
A left-invariant metric on $SO(3)$ with eigenvalues $\{\lambda_1^2,\lambda_2^2,1\}$ has nonnegative sectional curvature if and only if the following three inequalities hold:
\begin{align}\label{lamm}
2\lam_1^2 + 2\lam_2^2 - 3 + \lam_1^4 - 2\lam_1^2\lam_2^2 + \lam_2^4 &\geq 0,\notag\\
2\lam_1^2 - 3\lam_1^4 + 2\lam_1^2\lam_2^2 - 2\lam_2^2 + 1 + \lam_2^4 &\geq 0,\\
2\lam_2^2 - 3\lam_2^4 + 2\lam_1^2\lam_2^2 - 2\lam_1^2 + 1 + \lam_1^4 &\geq 0.\notag
\end{align}
\end{prop}
\begin{figure}[h!]
  \scalebox{0.3}{\includegraphics[angle=270]{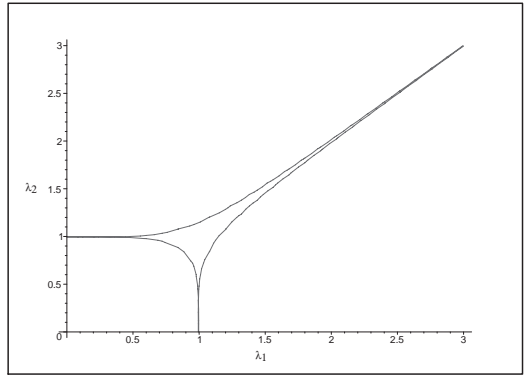}}
  \caption{The values of $\lam_1$ and $\lam_2$ for which $SO(3)$ has nonnegative sectional
          curvature, when $\lam_3=1$}
  \end{figure}
The intersection of the graph in Figure 1 with the identity function $\lambda_1=\lambda_2$ is interesting.  The eigenvalues $\{\lambda^2,\lambda^2,1\}$ yield nonnegative curvature if and only if $\lambda^2\geq 3/4$.  Re-scaled, this means that the eigenvalues $\{1,1,\mu^2\}$ yield nonnegative curvature if and only if $\mu^2\leq 4/3$.  Metrics with two equal eigenvalues on $SO(3)$ (or on its double-cover $S^3$) are commonly called \emph{Berger metrics}.  They are obtained by scaling the Hopf fibers by a factor of $\mu^2$.  The nonnegative curvature cut-off $\mu^2\leq 4/3$ is well-known.

\begin{proof}
Let $g_l$ be a left-invariant metric on $SO(3)$.  Let $\{u_1,u_2,u_3\}\subset\mg$ denote an oriented $g_0$-orthonormal basis of eigenvectors of $g_l$, with eigenvalues $\{\lambda_1^2,\lambda_2^2,\lambda_3^2\}$.  Then, 
$$\left\{e_1=\frac{u_1}{\lambda_1},e_2=\frac{u_2}{\lambda_2},e_3=\frac{u_3}{\lambda_3}\right\}$$ is a $g_l$-orthonormal basis of $\mg$.  The Lie bracket structure of $\mg$ is given by:
\begin{equation}\label{brack}[u_1,u_2]=u_3,\,\,[u_2,u_3]=u_1,\,\,[u_3,u_1]=u_2.\end{equation}
The structure constants, $\alpha_{ijk}=\lb[e_i,e_j],e_k\rb_l$, are:
$$\alpha_{123}=-\alpha_{213} = \frac{\lambda_3}{\lambda_1\lambda_2},\,\,\,\, \alpha_{231}=-\alpha_{321} =  \frac{\lambda_1}{\lambda_2\lambda_3},\,\,\,\, \alpha_{312}=-\alpha_{132} = \frac{\lambda_2}{\lambda_1\lambda_3},$$
with all others equal zero.  Equation~\ref{miln} gives:
\begin{multline*}
\kappa(e_1,e_2) = \frac{1}{2}\frac{\lam_3}{\lam_1\lam_2}\left(-\frac{\lam_3}{\lam_1\lam_2}+\frac{\lam_1}{\lam_2\lam_3}
 +\frac{\lam_2}{\lam_1\lam_3}\right) \\
 -\frac{1}{4}\left(\frac{\lam_3}{\lam_1\lam_2}-\frac{\lam_1}{\lam_2\lam_3}
  +\frac{\lam_2}{\lam_1\lam_3}\right)\left(\frac{\lam_3}{\lam_1\lam_2}
  +\frac{\lam_1}{\lam_2\lam_3}-\frac{\lam_2}{\lam_1\lam_3}\right).
\end{multline*}
When $\lambda_3=1$, the inequality $\kappa(e_1,e_2)\geq 0$ is equivalent to the first inequality of~\ref{lamm}.  The second and third inequalities are equivalent to $\kappa(e_2,e_3)\geq 0$ and $\kappa(e_1,e_3)\geq 0$ respectively.

It remains to prove the following: if the planes spanned by pairs of $\{e_1,e_2,e_3\}$ are nonnegatively curved, then all planes are nonnegatively curved.  It is immediate from P\"uttmann's formula~\ref{put} that $R(e_i,e_j)e_k=0$ whenever $i,j,k$ are distinct.  This implies that the curvature operator is diagonal in the basis $\{e_1\wedge e_2,e_2\wedge e_3,e_3\wedge e_1\}$.  The result follows.
\end{proof}

To classify up to isometry the left-invariant metrics of nonnegative sectional curvature on $SO(3)$ it remains to observe that:
\begin{prop} \label{HI}Two left-invariant metrics on $SO(3)$ are isometric if and only if they have the same eigenvalues.
\end{prop}
\begin{proof}
Suppose first that $f:(SO(3),g_1)\ra(SO(3),g_2)$ is an isometry, where $g_1$ and $g_2$ are left-invariant metrics.  We lose no generality in assuming that $f(I)=I$.  The linear map $\wedge^2(\mg)\ra\wedge^2(\mg)$ induced by $df_I:\mg\ra\mg$ must send eigenvectors of the curvature operator of $g_1$ to eigenvectors of the curvature operator of $g_2$.  Since these eigenvectors are wedge products of eigenvectors of the metrics, $df_I:\mg\ra\mg$ must send the eigenvectors of $g_1$ to the eigenvectors of $g_2$.  Since $df_I$ preserves lengths, the two metrics must have the same eigenvalues.

Conversely, suppose that $g_1$ and $g_2$ are left-invariant metric on $SO(3)$ with eigenvectors $E=\{e_1,e_2,e_3\}$ and $F=\{f_1,f_2,f_3\}$ respectively.  Assume that the eigenvalues of the two metrics agree.  The bases $E$ and $F$ can be chosen to be $g_0$-orthonormal and to have the same orientation.  There exists $x\in SO(3)$ such that $\Ad_x:\mg\ra\mg$ sends $E$ to $F$.  It is straightforward to verify that $C_x:(SO(3),g_1)\ra(SO(3),g_2)$ (conjugation by $x$) is an isometry.
\end{proof}

We conclude our investigation of $SO(3)$ by observing a relationship between the scalar curvature, $\rho$, of the metric with eigenvalues $\lambda_1^2,\lambda_2^2,\lambda_3^2$ and the area of a triangle with these three side-lengths.

Using notation from the proof of Proposition~\ref{P:rocket}, the scalar curvature is computed as
\begin{equation*}
\rho=2(\kappa(e_1,e_2)+\kappa(e_1,e_3)+\kappa(e_2,e_3)).
\end{equation*}
Summing our previous formulas for the curvatures of these planes yields:
\begin{equation*}
\rho=\frac{(-\lam_1+\lam_2+\lam_3)(\lam_1-\lam_2+\lam_3)(\lam_1+\lam_2-\lam_3)(\lam_1+\lam_2+\lam_3)}{2\lam_1^2\lam_2^2\lam_3^2}
\end{equation*}
Recall Heron's formula for the area, $A$, of a triangle with sides of length $a,b,c$:
\begin{equation*}
A=\sqrt{s(s-a)(s-b)(s-c)},
\end{equation*}
where $s=(a+b+c)/2$ is the semiperimeter of the triangle.  So when $\lam_1,\lam_2,\lam_3$ are possible lengths of a triangle, we have:
\begin{equation*}
\rho=\frac{8A^2}{\lam_1^2\lam_2^2\lam_3^2}.
\end{equation*}

\section{Metrics on $U(2)$}
Let $G=U(1)\times SU(2)$.  In this section, we classify the left-invariant metrics on $G$ with nonnegative curvature.  This is equivalent to solving the problem for $U(2)$, which has the same Lie algebra.

Let $g_0$ denote a bi-invariant metric on $G$.  Suppose that $g_l$ is a left-invariant metric on $G$.  The Lie algebra of $G$ is $\mg=u(1)\oplus su(2)$.  The factors $u(1)$ and $su(2)$ are orthogonal with respect to $g_0$. If they are orthogonal with respect to $g_l$, then $g_l$ is a product metric, so the problem reduces to the one from the previous section.

We will see, however, that nonnegatively curved metrics need not be product metrics.  This might be surprising, since when $g_l$ is nonnegatively curved, the splitting theorem (chapter 4 of~\cite{CG}) implies that the pull-back of $g_l$ to the universal cover $\R\times SU(2)$ is isometric to a product metric.  The subtlety is that the metric product structure and the group product structure need not agree.  Further, the metric $g_l$ on $G$ need only be locally isometric to a product metric, not globally.

Let $\{u_1,u_2,u_3\}$ denote an orthonormal basis of eigenvectors of the restriction of $g_l$ to the $su(2)$ factor.  Let $\{\lambda_1^2,\lambda_2^2,\lambda_3^2\}$ denote the corresponding eigenvalues, which we call ``restricted eigenvalues'' of $g_l$.  Notice that $$\left\{e_1=\frac{u_1}{\lambda_1},e_2=\frac{u_2}{\lambda_2},e_3=\frac{u_3}{\lambda_3}\right\}$$ is $g_l$-orthonormal.  Let $e_0\in\mg$ span the $u(1)$ factor.  Let $E_0$ span the $g_l$-orthogonal compliment of the $su(2)$ factor, and be $g_l$-unit length.

\begin{prop}\label{U2} $(G,g_l)$ has nonnegative curvature if and only if $\{\lambda_1,\lambda_2\,\lambda_3\}$ satisfy the restrictions of Proposition~\ref{P:rocket}, and one of the following conditions holds:
\begin{enumerate}
\item $E_0$ is parallel to $e_0$ (in this case, $g_l$ is a product metric)
\item $\lambda_1=\lambda_2=\lambda_3$ (in this case, $E_0$ is arbitrary)
\item $\lambda_1=\lambda_2$ and $E_0\in\text{span}\{e_0,e_3\}$ (or analogously if a different pair of $\lambda$'s agree).
\end{enumerate}
\end{prop}

We prove the proposition with a sequence of lemmas.

\begin{lem}\label{L2}
The map $ad_{E_0}:\mg\ra\mg$ is skew-adjoint with respect to $g_l$ if and only if one of the three conditions of Proposition~\ref{U2} holds.\end{lem}
\begin{proof}
Write $E_0=ae_0+be_1+ce_2+de_3$.  In the $g_l$-orthonormal basis $\{E_0,e_1,e_2,e_3\}$ of $\mg$,
the structure constants $\alpha_{0jk}=\lb[E_0,e_j],e_k\rb_l$ are as follows:
$$
\alpha_{012}=d\frac{\lambda_2}{\lambda_1},\,\,
\alpha_{021}=-d\frac{\lambda_1}{\lambda_2},\,\,
\alpha_{023}=b\frac{\lambda_3}{\lambda_2},\,\,
\alpha_{032}=-b\frac{\lambda_2}{\lambda_3},\,\,
\alpha_{031}=c\frac{\lambda_1}{\lambda_3},\,\,
\alpha_{013}=-c\frac{\lambda_3}{\lambda_1}.
$$
Here we assume for convenience that $g_0$ is scaled such that the Lie bracket structure of the $SU(2)$ factor is given by equation~\ref{brack}.  The lemma follows by inspection, since $ad_{E_0}$ is skew-adjoint if and only if $\alpha_{0ij}=-\alpha_{0ji}$ for all $i,j$.
\end{proof}
The next lemma applies more generally to any left-invariant metric on any Lie group.
\begin{lem}\label{par}The left-invariant vector field $x\in\mathfrak{g}$ is parallel if and only if $ad_x$ is skew-adjoint and $\langle x,\lbrack \mathfrak{g},\mathfrak{g} \rbrack\rangle_l=0$.
\end{lem}
\begin{proof}
For all $y,z\in\mg$,
\begin{equation}\label{cov}
2\langle y,\nabla_z x\rangle_l=\langle z,\lbrack y,x\rbrack\rangle_l+\langle x,\lbrack y,z\rbrack\rangle_l-\langle y,\lbrack x,z\rbrack\rangle_l.
\end{equation}

If $ad_x$ is skew-adjoint, then first and last terms of~\ref{cov} sum to $0$.  If additionally $\langle x,\lbrack \mathfrak{g},\mathfrak{g} \rbrack\rangle_l=0$, then the middle term is also $0$.  So $\nabla_zx=0$ for all $z\in\mg$, which means that $x$ is parallel.

Conversely, assume that $x$ is parallel, so the left side of~\ref{cov} equals $0$ for all $y,z\in\mg$.  When $y=z$, this yields $2\lb y,[y,x]\rb_l=0$ for all $y\in\mg$, which implies that $ad_x$ is skew-adjoint.  This property makes the first and third terms of~\ref{cov} sum to zero, so $\lb x,[y,z]\rb_l=0$ for all $y,z\in\mg$.  In other words, $\lb x,[\mg,\mg]\rb_l=0$.
\end{proof}

The next lemma from~\cite{milnor} also applies more generally to any left-invariant metric on any Lie group.  We use $r$ to denote the Ricci curvature of $g_l$.
\begin{lem}[Milnor]\label{L:mil}If $x\in\mg$ is $g_l$-orthogonal to the commutator ideal $[\mathfrak{g},\mathfrak{g}]$, then $r(x)\leq 0$, with equality if and only if $ad_x$ is skew-adjoint with respect to $g_l$.
\end{lem}

\begin{proof}[Proof of Proposition~\ref{U2}]
Suppose that $g_l$ has nonnegative curvature.  Then $ad_{E_0}$ is skew-adjoint by Milnor's Lemma~\ref{L:mil}.  Next, Lemma~\ref{L2} implies that one of the three conditions of the proposition hold.  It remains to prove that $\{\lambda_1,\lambda_2,\lambda_3\}$ satisfy the constraints for eigenvalues of a nonnegatively curved metric on $SO(3)$ (or equivalently on $SU(2)$).  By Lemma~\ref{par}, $E_0$ is parallel.  Using Equation~\ref{cov}, this implies that the $SU(2)$ factor of $G$ is totally geodesic, so its induced metric has nonnegative curvature, which gives the constraints on the $\lambda's$.

Conversely, suppose $\{\lambda_1,\lambda_2,\lambda_3\}$ satisfies the constraints for eigenvalues of a nonnegatively curved metric on $SU(2)$, and that one of the three conditions of the proposition holds.  By Lemma~\ref{L2}, $ad_{E_0}$ is skew-adjoint, so by Lemma~\ref{par}, $E_0$ is parallel.  This implies that the $SU(2)$ factor of $G$ is totally geodesic.  It has nonnegative curvature because of the constraints on the $\lambda$'s.  Consider the curvature operator of $g_l$ expressed in the following basis of $\wedge^2(\mg)$:
$$\{E_0\wedge e_1,E_0\wedge e_2,E_0\wedge e_3,e_1\wedge e_2,e_2\wedge e_3,e_3\wedge e_1\}.$$
Since $E_0$ is parallel, $R(E_0\wedge e_i)=0$.  Furthermore, $R(e_i\wedge e_j)$ is calculated in the totally geodesic $SU(2)$.  It follows that the curvature operator is nonnegative.
\end{proof}

We refer to metrics of types (2) and (3) in proposition~\ref{U2} as ``twisted metrics.''  We know from the splitting theorem that nonnegatively curved twisted metrics are locally isometric to untwisted (product) metrics.   We end this section by explicitly exhibiting the local isometry between a twisted metric and a product metric.

Suppose that $g_l$ is a twisted metric.  Recall that $\{E_0,e_1,e_2,e_3\}$ is a $g_l$-orthonormal basis of $\mg$.  Let $\tilde{g}_l$ denote the pull-back of $g_l$ to the universal cover, $\R\times SU(2)$, of $G$.  Let $g'$ denote the product metric on $\R\times SU(2)$ for which $\{e_0,e_1,e_2,e_3\}$ forms a $g'$-orthonormal basis of $\mg$.  Define
$$f:(\R\times SU(2),g')\ra(\R\times SU(2),\tilde{g}_l)$$
as follows:
\begin{eqnarray*}
f(t,g) & = & \text{flow from $(0,g)$ for time $t$ along $E_0$}\\
       & = & R_{e^{tE_0}}(0,g) = (at,ge^{t\hat{E}_0}),
\end{eqnarray*}
where $E_0 = ae_0 + be_1 + ce_2 + de_3$ and $\hat{E}_0 = be_1+ce_2+de_3$.
\begin{prop}\label{isom}
If $g_l$ has nonnegative curvature, then $f$ is an isometry.
\end{prop}
\begin{proof} First,
\begin{equation}\label{df1} df_{(t,g)}(e_0) = \frac{d}{ds}\Big|_{s=0}\left(a(t+s),ge^{(t+s)\hat{E}_0}\right) = \left(a,\left(ge^{t\hat{E}_0}\right)\hat{E}_0\right),\end{equation}
which is the value at $f(t,g)$ of the left-invariant vector field $E_0$.  So $df$ send the left-invariant field $e_0$ to the left-invariant field $E_0$.  Next, for $i\in\{1,2,3\}$,
\begin{equation}
\label{df2}df_{(t,g)}(e_i g) = \frac{d}{ds}\Big|_{s=0} f\left(t,e^{se_i}g\right) =  \frac{d}{ds}\Big|_{s=0}\left(at,e^{se_i}ge^{t\hat{E}_0}\right)
=\left(0,e_i\left(ge^{t\hat{E}_0}\right)\right),\end{equation}
which is the value at $f(g,t)$ of the right-invariant vector field determined by $e_i$.  In other words, $df$ sends the right-invariant vector field determined by $e_i$ to itself.

In summary,~\ref{df1} implies that $df$ sends the $g'$-unit-length field $e_0$ to the $\tilde{g}_l$-unit-length field $E_0$.  Further,~\ref{df2} implies that $df$ sends the $g'$-orthogonal compliment $e_0$ to the $\tilde{g}_l$-orthogonal compliment of $E_0$ (both of which equal $\text{span}\{e_1,e_2,e_3\}$).  It remains to verify that the restriction $df:\text{span}\{e_1,e_2,e_3\}\ra\text{span}\{e_1,e_2,e_3\}$ is an isometry.

Under condition (2) of Proposition~\ref{U2}, $g'$ restricts to a bi-invariant metric on the $SU(2)$ factor.  Under condition (3), $g'$ restricts to a left-invariant and $\Ad_H$-invariant metric, where $H=\exp(\text{span}(\hat{E}_0))$.  In either case, the metric has exactly enough right-invariance to give the desired result from Equation~\ref{df2}.
\end{proof}
Proposition~\ref{HI} and~\ref{isom} together imply that:
\begin{cor}
Two left-invariant metrics with nonnegative curvature on $U(1)\times SU(2)$ are locally isometric if and only if their restricted eigenvalues $\{\lambda_1,\lambda_2,\lambda_3\}$ are the same.
\end{cor}
These local isometries are really isometries between the universal covers.  In general, they are not group isomorphisms. They do not generally descend to global isometries between twisted and product metrics on $U(1)\times SU(2)$.
\section{Review of Cheeger's method}
In the literature, the only examples of left-invariant metrics with nonnegative curvature on compact Lie groups come from a construction which, in its greatest generality, is due to Cheeger~\cite{Cheeger}.  In this section, we review Cheeger's method.

Let $G$ be a compact Lie group, and let $H\subset G$ be a closed subgroup.  Let $\mh\subset\mg$ denote their Lie algebras.  Let $g_0$ be a left-invariant and $\Ad_H$-invariant metric on $G$ with nonnegative curvature.  Let $g_H$ be a right-invariant metric on $H$ with nonnegative curvature.  Denote:
$$\Delta H:=\{(h,h)\in H\times G\mid h\in H\}.$$
The right-action of $\Delta H$ on $(G,g_0)\times (H,g_H)$ is by isometries.  The quotient, $(G\times H)/\Delta H$, is diffeomorphic to $G$ via the diffeomorphism $[g,h]\stackrel{f}{\mapsto} g\cdot h^{-1}$.  This quotient inherits a Riemannian submersion metric, $g_1$, with nonnegative curvature.  We write:
$$(G,g_1) = ((G,g_0)\times (H,g_H))/\Delta H.$$
In fact, $g_1$ is a left-invariant metric on $G$, since for all $a\in G$, $f^{-1}$ associates $L_a:G\ra G$ (left-multiplication by $a$) to the isometry $[g,h]\mapsto [ag,h]$.  Similarly, if $g_H$ is bi-invariant, then $g_1$ is $\Ad_H$-invariant, since for all $a\in H$, $f^{-1}$ associates $R_a:G\ra G$ to the isometry $[g,h]\mapsto [g,a^{-1}h]$.

The metrics $g_0$ and $g_1$ agree orthogonal to $\mh$.  In other words, If $X\in\mg$ is $g_0$-orthogonal to $\mh$, then it is $g_1$-orthogonal to $\mh$, and $g_0(X,X)=g_1(X,X)$.  It remains to describe $g_1$ on $\mh$.  For this, let $\{e_1,...,e_k\}$ denote a $g_H$-orthonormal basis of $\mh$, and regard the $e$'s as vectors in $\mg$.  Let $A$ and $\tilde{A}$ respectively denote the $k$-by-$k$ matrices whose entries are:
$$a_{ij} = g_0(e_i,e_j),\,\,\,\,\,\tilde{a}_{ij} = g_1(e_i,e_j).$$
The restrictions of $g_1$ and $g_0$ to $\mh$ are related by the equation:
\begin{equation}\label{Cheg}\tilde{A} = A(I+A)^{-1}.\end{equation}

A common special case occurs when $g_0$ is bi-invariant, and $g_H$ is a multiple, $\lambda$, of the restriction of $g_0$ to $H$.  In this case, which is studied in~\cite{Eschenburg}, $g_1$ can be described as follows:
\begin{equation}\label{esh}g_1(X,Y) = g_0(X^{\mp},Y^{\mp}) + t\cdot g_0(X^{\mh},Y^{\mh}),\end{equation}
where $X^{\mp}$ (respectively $X^\mh$) denotes the $g_0$-projection of $X$ orthogonal to (respectively onto) $\mh$, and
$$t=\frac{\lambda}{\lambda+1}.$$
Thus, $g_1$ is obtained from $g_0$ by uniformly shrinking all vectors in $\mh$.

In this case, $g_1$ is $\Ad_H$-invariant, and is therefore $\Ad_K$-invariant for any closed subgroup $K\subset H$.  So Cheeger's method can be applied again:
$$(G,g_2) = ((G,g_1)\times(K,g_K))/\Delta K.$$
If $g_K$ is a multiple of the restriction of $g_0$ to $K$, then the process can be again repeated for any subgroup of $K$, and so on.

In summary, whenever $H_0\subset H_1\subset \cdots\subset H_{l-1}\subset H_l=G$ is a chain of closed subgroups of $G$ with Lie algebras $\mh_0\subset\mh_1\subset\cdots\subset\mh_{l-1}\subset\mh_l=\mg$, one can apply Cheeger's method $l$-times.  One chooses a starting bi-invariant metric $g_0$ on $G$ and $l$ constants $\{\lambda_0,...,\lambda_{l-1}\}$.  The result is a new left-invariant metric with nonnegative curvature on $G$.  The sub-algebras $\mh_i$ are separately scaled by factors determined by the $\lambda$'s, with $\mh_i$ scaled a greater amount than $\mh_{i+1}$.  More precisely, the eigenvalues of the metric are a strictly increasing sequence $$t_0<t_1<\cdots<t_{l-1}<t_l=1.$$  The $\lambda's$ can be chosen to give any such strictly increasing sequence.  The eigenspace of $t_0$ equals $\mh_0$.  For $i>1$, the eigenspace of $t_i$ equals the $g_0$-orthogonal compliment of $\mh_{i-1}$ in $\mh_i$.

\section{Do all metrics come from Cheeger's method?}
Let $G$ be a compact Lie group.  Does every left-invariant metric on $G$ with nonnegative curvature come from Cheeger's method?  In this section, we outline several ways to precisely formulate this question, and address the cases $G=SO(3)$ and $G=U(2)$.

First, one might ask whether all examples arise by starting with a bi-invariant metric $g_0$ on $G$ and applying Cheeger's method to a chain of subgroups $H_0\subset\cdots\subset H_{l-1}\subset G$, each time choosing the metric on $H_i$ to be a multiple, $\lambda_i$, of the restriction of $g_0$ to $H_i$.  The answer is clearly no, since not all metrics on $SO(3)$ arise in this fashion.  But $SO(3)$ is a very special case, since its bi-invariant metric has positive curvature.  There are several ways to modify the question to reflect the guess that $SO(3)$ is the only exception.

What if, in each application of Cheeger's method, one allows more general right-invariant metrics on the $H_i$'s?   For example, if $H_i\cong SO(3)$ or $SU(2)$, we allow any right-invariant metric with nonnegative curvature.  With this added generality, the answer is still no, since metrics on $SO(3)$ which are nonnegatively but not positively curved do not arise in this fashion:
\begin{prop} If $g_o$ is a bi-invariant metric on $SO(3)$ and $g_R$ is a right-invariant metric with nonnegative curvature on $SO(3)$, then the following has strictly positive curvature:
$$(SO(3), g_l):=((SO(3),g_0)\times(SO(3),g_R))/\Delta SO(3).$$
Further, every positively-curved left-invariant metric $g_l$ can be described in this way for some bi-invariant metric $g_0$ and some positively curved right-invariant metric $g_R$.
\end{prop}
\begin{proof}
Let $\pi:(SO(3),g_0) \times (SO(3), g_R)\ra(SO(3), g_l)$ denote the projection, which is a Riemannian submersion.  Let $\{e_1,e_2,e_3\}\subset so(3)$ be a $g_0$-orthonormal basis of eigenvectors of the metric $g_R$, with eigenvalues denoted $\{\lambda_1,\lambda_2,\lambda_3\}$.  The horizontal space of $\pi$ at the identity is the following subspace of $so(3)\times so(3)$:
\begin{align}\label{wu1}
  \Hor_{(e,e)} &= \{ (V,W)\mid \lb(V,W), (e_i, e_i)\rb=0 \text{ for all }i=1,2,3\} \\
       &= \{ (V,W)\mid \lb V,e_i\rb_0 +\lambda_i \lb W,e_i\rb_0=0 \text{ for all }i=1,2,3\}. \notag\end{align} 
Now let $X_1=(V_1, W_1), X_2=(V_2, W_2) \in \Hor_{(e,e)}$ be linearly independent vectors.  From equation \ref{wu1}, we see that $V_1, V_2$ are linearly independent as well.  Since $(SO(3), g_0)$ has positive curvature,  $\kappa(V_1, V_2)>0$, so $\kappa(X_1, X_2)>0$.  Finally, O' Neill's formula implies that $(SO(3), g_l)$ has strictly positive curvature. 

To prove the second statement of the proposition, notice that by equation~\ref{esh}, the eigenvalues $\{\tilde{\lambda}_1,\tilde{\lambda}_2,\tilde{\lambda}_3\}$ of the metric $g_l$ are determined from the eigenvalues $\{\lambda_1,\lambda_2,\lambda_3\}$ of the metric $g_R$ as follows:
\begin{equation}\label{Matt}\tilde{\lambda}_i=\frac{\lambda_i}{1+\lambda_i}\text{ for each }i=1,2,3.\end{equation}
Suppose that $(\tilde{\lambda}_1,\tilde{\lambda}_2)$ is an arbitrary pair such that the triplet $\{\tilde{\lambda}_1,\tilde{\lambda}_2,1\}$ strictly satisfies the inequalities of Proposition~\ref{P:rocket}.  Define:
$$\lambda_1:=\frac{a\tilde{\lambda}_1}{1+a(1-\tilde{\lambda}_1)},\,\,\,
\lambda_2:=\frac{a\tilde{\lambda}_2}{1+a(1-\tilde{\lambda}_2)},\,\,\,\lambda_3:=a.$$
For small enough $a$, the triplet $\{\lambda_1,\lambda_2,\lambda_3\}$ strictly satisfies the inequalities of Proposition~\ref{P:rocket}.  This is because, when scaled such that the third equals $1$, the triplet approaches $\{\tilde{\lambda}_1,\tilde{\lambda}_2,1\}$ as $a\ra 0$.
If $g_R$ has eigenvalues $\{\lambda_1,\lambda_2,\lambda_3\}$, then by equation~\ref{Matt}, $g_l$ has the following eigenvalues:
$$\left\{\left(\frac{a}{1+a}\right)\tilde{\lambda}_1, \left(\frac{a}{1+a}\right)\tilde{\lambda}_2,\left(\frac{a}{1+a}\right)\right\}.$$
Further, $g_R$ has positive curvature because the inner products which extend to right-invariant metrics with positive curvature are the same ones that extend to left-invariant metrics of positive curvature.  This is because $a\mapsto a^{-1}$ is an isometry between left and right-invariant metrics determined by the same inner product.

Thus, $g_R$ can be chosen such that $g_l$ is a multiple of any prescribed positively curved metric.  Then, by scaling the bi-invariant metric $g_0$, this multiple can be made to be $1$.
\end{proof}

Starting with a bi-invariant metric $g_0$ on $G=SU(2)\times U(1)$ and applying Cheeger's method, one can not obtain product metrics for which the $SU(2)$-factor has nonnegative but not positive curvature, nor can one obtain twisted metrics of type (2).  Notice that the only chains of increasing-dimension subgroups are $U(1)\subset G$ (possibly embedded diagonally), $T^2\subset G$ (any maximal torus), $U(1)\subset T^2\subset G$, and $SU(2)\subset G$.

To obtain more metrics, in addition to allowing general right-invariant metrics on the $H_i$'s as above, one could allow a more general starting metric $g_0$.  For example, when $G=SU(2)\times U(1)$, one could allow $g_0$ to be any left-invariant product metric with nonnegative curvature.  Even with this added generality, it is straightforward to see that one cannot obtain twisted metrics of type (3) for which the totally geodesic $SU(2)$ has nonnegative but not positive curvature.

\bibliographystyle{amsplain}

\end{document}